\documentclass{mcom-l}
\usepackage{url}
\usepackage{booktabs}
\newcommand{\Fp}{\mathbb{F}_p}
\newcommand{\Fpp}{\mathbb{F}_{p'}}
\newcommand{\Fptwo}{\mathbb{F}_{p^2}}
\newcommand{\Fpt}{\mathbb{F}_{p^3}}
\newcommand{\Fpn}{\mathbb{F}_{p^n}}
\newcommand{\Fq}{\mathbb{F}_q}
\newcommand{\Fqbar}{\overline{\mathbb{F}}_q}
\newcommand{\Fqk}{\mathbb{F}_{q^k}}

\newcommand{\Fqd}{\mathbb{F}_{q^d}}
\newcommand{\Fqtwo}{\mathbb{F}_{q^2}}

\newcommand{\orda}{|\alpha|}

\newcommand{\identity}{1_{\scriptscriptstyle{G}}}

\newcommand{\Pic}{{\rm Pic^0}}

\newcommand{\Ctwist}{\tilde{C}}

\newtheorem{theorem}{Theorem}
\newtheorem{proposition}{Proposition}

\newtheorem{lemma}{Lemma}
\newtheorem{algorithm}{Algorithm}

\theoremstyle{definition}
\newtheorem{definition}{Definition}

\theoremstyle{remark}

\numberwithin{equation}{section}

\begin{document}

% \title[short text for running head]{full title}
\title{A generic approach to searching for Jacobians}

%    Only \author and \address are required; other information is
%    optional.  Remove any unused author tags.

%    author one information
% \author[short version for running head]{name for top of paper}
\author{Andrew V. Sutherland}
\address{Massachusetts Institute of Technology}
\email{drew@math.mit.edu}

%    \subjclass is required.

%    Abstract is required.
\begin{abstract}
We consider the problem of finding cryptographically suitable Jacobians.  By applying a probabilistic generic algorithm to compute the zeta functions of low genus curves drawn from an arbitrary family, we can search for Jacobians containing a large subgroup of prime order.  For a suitable distribution of curves, the complexity is subexponential in genus 2, and $O(N^{1/12})$ in genus 3.  We give examples of genus 2 and genus 3 hyperelliptic curves over prime fields with group orders over $180$ bits in size, improving previous results.  Our approach is particularly effective over low-degree extension fields, where in genus 2 we find Jacobians over $\Fptwo$ and trace zero varieties over $\Fpt$ with near-prime orders up to 372 bits in size.  For $p = 2^{61}-1$, the average time to find a group with 244-bit near-prime order is under an hour on a PC.
\end{abstract}

\maketitle

\section{Introduction}

Algebraic curves over finite fields have proven to be a fertile source of groups for cryptographic applications based on the discrete logarithm problem.  This has spurred the development of highly efficient algorithms for group computation that are now available for many types of curves, including hyperelliptic, superelliptic, Picard, and $C_{ab}$ curves \cite{Basiri:C34Arithmetic,Basiri:Superelliptic,Flon:PicardArithmetic,Flon:NonHyperGenus3,Gaudry:FastGenus2Arithmetic,Khuri:GeneralCurves,Lange:Genus2Arithmetic,Wollinger:Genus3Arithmetic}.
The group of interest consists of the $\Fq$-rational points on the Jacobian variety of a curve $C$, or, equivalently, the divisor class group of degree 0, $\Pic(C)$.  We denote this group $J(C/\Fq)$, or simply $J(C)$, and call it the \emph{Jacobian} of C.

For cryptographic use one typically seeks Jacobians with near-prime orders in the range $2^{160}$ to $2^{256}$, and we also consider subgroups of Jacobians that offer comparable (perhaps superior) performance and security parameters, such as trace zero varieties \cite{Cohen:HECHECC,Frey:TraceZero,Lange:TraceZero}.  The existence of various index calculus algorithms has centered attention on hyperelliptic curves of genus $g \le 3$ \cite{Diem:PlaneCurveIndexCalculus,Enge:HighGenusIndexCalculus,Gaudry:SmallGenusIndexCalculus}.  We similarly focus on the hyperelliptic case, although our results may be applied to any family of low genus curves.  

To assess the cryptographic suitability of a group, it is necessary to know its order.  For curves of genus 1 (elliptic curves), several effective point-counting algorithms are available.  The most general are $\ell$-adic methods, based on Schoof's algorithm \cite{Schoof:ECPointCounting1,Schoof:ECPointCounting2}, and for small characteristic fields there are more efficient $p$-adic methods \cite{Kedlaya:PointCounting,Kedlaya:ZetaFunctions,Mestre:AGMPointCounting,Satoh:ECPointCounting}.  The $p$-adic methods, particularly Kedlaya's algorithm \cite{Kedlaya:PointCounting}, readily extend to higher genus curves, and have proven effective over small and medium characteristic fields \cite{Gaudry:PointCountingMediumCharacterstic,Vercauteren:ZetaFunctions}.  Generalization of the $\ell$-adic methods has been more difficult.  The best results for genus 2 curves over prime fields report roughly a week to compute the order of a group of size $\approx 2^{164}$ \cite{Gaudry:Genus2PointCounting}.  In genus 3, no effective $\ell$-adic methods are available, however, recent work on extending $p$-adic methods to larger characteristic fields has enabled computation of group orders up to size $\approx 2^{150}$ over prime fields \cite{Bostan:LinearRecurrences,Harvey:LargeCharacteristicKedlaya}.  In both cases the algorithms are memory intensive, limiting their applicability to larger groups.  There are other methods for curves with special properties \cite{Furukawa:SpecialHECCount,Weng:GMGenus3,Weng:CMGenus2}, but for general curves of genus $g > 1$ over large characteristic fields, efficiently finding cryptographically suitable Jacobians remains an open problem \cite{Cohen:HECHECC,Lange:OpenProblems}.\footnote{There is a polynomial-time $\ell$-adic algorithm due to Pila for arbitrary abelian varieties  \cite{Pila:PointCounting}, but it is not practical for groups of cryptographic size.}
% \cite[p. 4]{Lange:OpenProblems}\cite[p. 568]{Cohen:HECHECC}.

The solution we propose is probabilistic.  As a point-counting algorithm, it usually fails.  However, given a suitably diverse family of curves, it succeeds often enough to be effective as a search algorithm, with (heuristically) sub-exponential performance.  One can then search for curves with a desired property, such as cryptographic suitability.

\section{Overview}

We apply neither $\ell$-adic nor $p$-adic methods, relying instead on \emph{generic} algorithms.  These perform group computations in a representation-independent manner, needing only a ``black box" to implement the group law.  As a result of the work cited above, we have many highly efficient black boxes at our disposal.

The $\Theta(\sqrt{N})$ complexity of birthday-paradox algorithms, such as Pollard's rho method \cite{Pollard:RhoDL} and Shanks' baby-steps giant-steps algorithm \cite{Shanks:BabyGiant}, makes them too slow to effectively compute the order of a large group, even when fairly tight bounds on the order are known (as in  \cite{Stein:HyperellipticBounds}).
Alternatively, one may apply a generic version of Pollard's $p-1$ technique \cite{Pollard:Pminus1}, exponentiating by many small primes.  This can be quite effective if the group order happens to be smooth (no large prime factors), but the worst case complexity is $\Theta(N)$.

Surprisingly, a combination of these two generic approaches is faster than either alone.  The author's thesis \cite{Sutherland:Thesis} presents a $o(\sqrt{N})$ algorithm to compute the order of an element in any finite group.  For a  family of abelian groups with orders uniformly distributed over a large interval, the average running time is $O(\sqrt{N}/\log{N})$.  This average is dominated by a $o(1)$ fraction of worst cases (groups of prime order).  The median complexity is $O(N^{0.344})$, and when the group order is highly composite, the algorithm is very fast.  By applying a bounded amount of computation to each group in a family, we can hope to find one whose order is easily computed.   This approach is similar to some algorithms for integer factorization \cite{Lenstra:ECM,Schnorr:Factorization}, and has been successfully applied to compute ideal class groups of imaginary quadratic number fields with 100-digit discriminants. For a suitable distribution of group orders, the complexity is subexponential: $L(1/2,\sqrt{2})= O\bigl(\exp\bigl[\bigl(\sqrt{2}+o(1)\bigr)(\log{N}\log\log{N})^{1/2}\bigr]\bigr)$.
%\footnote{$L(a,c) $}
% for time and $L(1/2,\sqrt{2}/2)$ for space  \cite[Proposition 5.5]{Sutherland:Thesis}.\footnote{The space complexity can be made polynomial, but this is rarely necessary.}

The ability to quickly find a group with highly composite order would seem little use in the search for cryptographically suitable groups, as this is precisely the opposite of what is desired.  However, the order of one Jacobian may be used to compute the orders of several others.  For low genus curves, given $\#J(C)$, we can readily recover the zeta function of $C$ (\ref{equation:ZetaDefinition}) via the $L$-polynomial $P(z)$ that appears in its numerator.  This process is trivial in genus 1, almost trivial in genus 2, and in genus 3 we give a generic algorithm requiring $O(N^{1/12})$ group operations.\footnote{This is asymptotically exponential, but negligible for the size groups we consider.}

With $P(z)$ in hand, we can compute the order of $J_d(C) = J(C/\Fqd)$ for any degree $d$ extension field (Lemma \ref{lemma:ZetaCount}).  When $b$ divides $a$, the abelian group $J_b(C)$ is a subgroup of $J_a(C)$, and we consider groups of the form
\begin{equation}
J_{a/b}(C) \cong J_a(C)/J_b(C).
\end{equation}
When $C$ is hyperelliptic, $J_{2/1}(C)$ is isomorphic to $J(\Ctwist)$, where $\Ctwist$ is the quadratic twist of the curve $C$ over $\Fq$ (Lemma \ref{lemma:Twist}).  In general, $J_{a/b}$ need not correspond to the Jacobian of a curve, however, we can compute in $J_{a/b}(C)$ using the group operation of $J_a(C)$.  The group $J_{a/1}(C)$ corresponds to a \emph{trace zero variety} (Lemma \ref{lemma:TraceZero}).  As noted by Lange, computation in the trace zero variety is typically more efficient than computation in a Jacobian of comparable size and genus, due to optimizations enabled by the Frobenius endomorphism \cite{Lange:TraceZero}.

The shape of the integer $\#J_{a/b}(C)$ is not particularly correlated with $\#J(C)$, and it is possible that the former may be near prime while the latter is quite smooth (Table \ref{table:Distribution}).  We should remark that this situation is not believed to diminish the security of the group $J_{a/b}(C)$.  Indeed, most of the elliptic curves in the NIST Digital Signature Standard \cite{NIST:FIPS186}, and nearly all the Certicom challenge curves \cite{Certicom:ECCChallenge}, have relatively insecure quadratic twists.\footnote{In particular, the order of the twist of the Certicom curve ECCp-163 contains no prime factors larger than $2^{56}$.  Discrete logarithms in the twist can be computed in well under an hour on a typical PC, yet the US\$30,000 cash prize remains unclaimed after nearly a decade.}

When considering the cryptographic use of Jacobians over extension fields, one must take into account the existence of transfers (cover attacks) which may reduce their effective security (Proposition \ref{proposition:Security}).  For hyperelliptic curves in genus 3 this essentially limits us to $J_{2/1}(C)\cong J(\Ctwist)$ over a prime field.
%For Picard curves, we may consider cubic twists, but note the attack of \cite{Diem:PlaneCurveIndexCalculus} on non-hyperelliptic curves.
For genus 2 curves, however, there are four groups with potentially competitive performance/security ratios when $q$ is prime:
\begin{equation}\label{equation:candidates}
 J_{2/1}(C)\cong J(\Ctwist),\qquad J_{3/1}(C),\qquad J_{3/1}(\Ctwist),\qquad J_{4/2}(C)\cong J(\Ctwist_2).
\end{equation}
Here $\Ctwist$ is the quadratic twist of $C$ in $\Fq$ and $\Ctwist_2$ is the quadratic twist in $\Fqtwo$.  The group $J(\Ctwist)$ has size $\approx q^2$, while the last three groups are each of size $\approx{q^4}$.
%Adjusted to reflect known attacks, the groups $J_{3/1}(C)$ and $J_{3/1}(\Ctwist)$ each offer security comparable to a genus 2 Jacobian over a prime field of size $O(q^{10/3})$ or $O(q^{16/5})$.  For $J_{2/4}(C)$ the comparable size is $O(q^3)$.  

To give a brief example, we used this approach on a family of $10^6$ random curves over $\Fp$, with $p = 2^{61}-1$ implying $\#J(C)\approx2^{122}$.  For suitably chosen parameters, a single PC (2.5 GHz AMD Athlon-64) test about two curves per second, successfully computing $\#J(C)$ for some curve in the family every four or five minutes on average.  We computed the zeta functions of some 2000 curves, finding 220 groups with 244-bit near-prime orders (cofactor $<5\%$ of the bits), including many of prime order.  For cryptographic use, most of these groups should be compared to genus 2 Jacobians over a prime field with a group size of 180 to 200 bits (Proposition \ref{proposition:Security}). Depending on the implementation, they may offer superior performance.

More detailed examples are provided in Section \ref{section:Examples}, including results for much larger groups.  We have successfully applied this approach over prime fields with $\#J(C)$ up to 186 bits in genus 2, and 183 bits in genus 3, improving previous bests of 164 bits (Gaudry and Schost \cite{Gaudry:Genus2PointCounting}) and 150 bits (Harvey \cite{Harvey:LargeCharacteristicKedlaya}).  We note, however, that while our method computes zeta functions of curves drawn from an arbitrary family, it will likely fail on any particular curve and should be distinguished from point-counting algorithms.  The new algorithm is faster, less memory intensive, and well suited to distributed implementation. 

\section{Mathematical Background}\label{section:Facts}
For a projective curve $C$ defined over $\Fq$, let $N_k$ count the points on $C$ in $\mathbb{P}(\Fqk)$.  The \emph{zeta function} of $C$ is the formal power series
\begin{equation}\label{equation:ZetaDefinition}
Z(C/\Fq,z) = \exp\left(\sum_{k=1}^{\infty}N_kz^k/k\right).
\end{equation}
Our interest in the zeta function stems from the well-known theorem of Weil \cite{Weil:ZetaFunction}, which we restrict here to projective curves defined over $\Fq$.  Henceforth we assume all curves are non-singular and irreducible over the algebraic closure $\Fqbar$.
\begin{theorem}[Weil]\label{theorem:Weil}
Let $C$ be a genus $g$ curve defined over $\Fq$.
\begin{enumerate}
\item
$Z(C/\Fq,z) = P(z)/\left[(1-z)(1-qz)\right],$ where $P(z) = \sum_{i=0}^{2g}a_iz^i$ has integer coefficients satisfying
$a_0 = 1$ and $a_{2g-i} = q^{g-i}a_i,$ for $0 \le i < g$.
\vspace{6pt}
\item
$P(z) = \prod_{i=1}^{2g}(1-\alpha_iz),\space \text{ with } |\alpha_i| = \sqrt{q}.$
\vspace{6pt}
\item
$N_k = q^k +1 - \sum_{i=1}^{2g}\alpha_i^k.$
\end{enumerate}
\end{theorem}
A proof can be found in chapters 8 and 10 of \cite{Lorenzini:ArithmeticGeometry}.  We call $P(z)$ the \emph{$L$-polynomial} of the curve $C$.  From (2) we obtain the bounds
\begin{equation}\label{equation:abounds}
|a_i| \le \binom{2g}{i}q^{i/2}.
\end{equation}
Let $J(C/\Fqk)$ denote the group of $\Fqk$-rational points on the Jacobian variety of $C$.  We write $J_k(C)$ for $J(C/\Fqk)$, and $J(C)$ for $J(C/\Fq)$, when $\Fq$ is understood.
\begin{lemma}\label{lemma:ZetaCount}
Let $C$ and $P(z)$ be as in Theorem \ref{theorem:Weil}.  Then
$$\#J_k(C) = \prod_{i=1}^k P(\omega^i),$$
where $\omega = e^{2\pi i/k}$ is a principal kth root of unity.
\end{lemma}
See \cite[8.5.12 and 8.6.2-3]{Lorenzini:ArithmeticGeometry} for a proof.  In particular, $\#J(C) = P(1)$, and applying (2) of Theorem 1 gives the \emph{Weil interval}:
\begin{equation}\label{equation:WeilInterval}
(\sqrt{q}-1)^{2g} \le \#J(C) \le (\sqrt{q}+1)^{2g}.
\end{equation}

The Frobenius automorphism $a\to a^q$ of the finite field $\Fq$ gives rise to a group endomorphism on $J_k(C)$, which we denote $\phi_q$.  The elements of $J_k(C)$ fixed by $\phi_q$ are precisely the subgroup $J(C)$.
When $d$ divides $k$, we define $J_{k/d}$ to be the image of $\phi_q^d-1$ on $J_k$, where 1 denotes the identity map.  Thus $J_{k/d}$ is a subgroup of $J_k(C)$ isomorphic to $J_k(C)/J_d(C)$.  We define the \emph{trace zero variety} of $J_k(C)$, denoted $T_k(C)$, to be the kernel of the group endomorphism
\begin{equation}
\phi_q^{k-1} +\phi_q^{k-2} +\cdots + 1.\notag
\end{equation}
\begin{lemma}\label{lemma:TraceZero}
Let $C$ be a curve defined over $\Fq$. Then $J_{k/1}(C) \subseteq T_k(C)$.  Equality holds when $J(C)$ is $k$-torsion free.
\end{lemma}
\begin{proof}
To show inclusion, we factor $\phi_q^k-1$ in the endomorphism ring of $J_k(C)$:
$$J_k(C) = \rm{ker}(\phi_{q^k}-1) = \rm{ker}(\phi_q^k-1) = \rm{ker}\left[(\phi_q-1)(\phi_q^{k-1} +\cdots + 1)\right].$$
The image of $\phi_q-1$ lies in the kernel of $\phi_q^{k-1} +\cdots + 1$, hence $J_{k/1}(C) \subseteq T_k(C)$.

When $J(C)$ is $k$-torsion free, $J(C)\cap T_k(C)$ is trivial and
$$\#J_k(C) = \#J(C)\#J_{k/1}(C) \le \#J(C)\#T_k(C) \le \#J_k(C),$$
implying equality.
\end{proof}

If $\Fq$ is a finite field of odd characteristic, we may define a \emph{hyperelliptic curve} of genus $g$ as a projective plane curve $C$ with affine part given by
\begin{equation}\notag
y^2 = f(x),
\end{equation}
where $f(x) \in \Fq[x]$ has non-zero discriminant and degree $d=2g+1$ or $2g+2$.\footnote{The black boxes we use assume $f$ is monic and $d=2g+1$ but our results do not require this.}  If $\alpha\in\Fq$ is a not a quadratic residue, the \emph{quadratic twist} of $C$ in $\Fq$, denoted $\Ctwist$, has affine part $\alpha y^2 = f(x)$, more conveniently expressed as
\begin{equation}
y^2 = \alpha^{d}f(x/\alpha).\notag
\end{equation}
Any choice of non-residue in $\Fq$ yields a curve isomorphic to $\Ctwist$.
\begin{lemma}\label{lemma:Twist}
Let $C$ be a hyperelliptic curve with $L$-polynomial $P(z)$.
\begin{center}
(1) $P(-z)$ is the $L$-polynomial of $\Ctwist$; $\qquad$ (2) $J_{2/1}(C) \cong J(\Ctwist)$.
\end{center}
\end{lemma}
\begin{proof}
Let $\tilde{N}_k$ count the points on $\tilde{C}$ in $\mathbb{P}(\Fqk)$.   For $k$ even, $\tilde{N}_k = N_k$, and for $k$ odd, $N_k+\tilde{N}_k = 2(q^k+1)$.  Applying (3) and (2) of Theorem \ref{theorem:Weil} proves (1), and Lemma \ref{lemma:ZetaCount} implies that $\#J_2(C) = P(1)P(-1) = \#J(C)\#J(\Ctwist)$.  The curves $C$ and $\Ctwist$ are isomorphic over $\Fqtwo$, thus we may regard $J(\Ctwist)$ as a subgroup of $J_2(C)$ that intersects $J(C)$ trivially.  $J_2(C)$ is then the (internal) product of the subgroups $J(C)$ and $J(\Ctwist)$, hence $J_{2/1}(C) \cong J_2(C)/J(C)\cong J(\Ctwist)$.
%thus $J_2(C)\cong J(C)\otimes J(\Ctwist)$ and $J_{2/1}(C) \cong J_2(C)/J(C)\cong J(\Ctwist)$.
\end{proof}

Finally, we note results that impact the effective security of hyperelliptic curves \cite{Diem:GHSAttack,Diem:CoverAttacks,Gaudry:GHSAttack,Gaudry:SmallGenusIndexCalculus,
Gaudry:AbelianVarietiesIndexCalculus,Smith:Genus3DiscreteLog,Theriault:IndexCalculus,Theriault:WeilDescent}. Taking a pessimistic view, we list the strongest potential attacks known at the time of writing.
\begin{proposition}\label{proposition:Security}
Let $C$ be a hyperelliptic curve of genus $g$ over $\Fq = \Fpn$.
\begin{enumerate}
\item
Discrete logarithms in $J(C)$ can be computed in time $O(p^{2-4/(2ng+1)})$. \cite{Gaudry:AbelianVarietiesIndexCalculus}
\item
If $J(C)$ is cyclic and $g\ge 3$ discrete logarithms in $J(C)$ can be computed in time $\tilde{O}(q^{2-2/g})$.  Heuristically, $J(C)$ need not be cyclic. \cite{Gaudry:SmallGenusIndexCalculus}
\item
If $C$ is a genus 2 curve, the discrete logarithm problem in the trace zero variety $T_3(C)$ may be transferred to the Jacobian of genus 6 curve over $\Fq$.  If $\#J_2(C)$ is divisible by 3, the genus may be 5. \cite{Diem:CoverAttacks}
%A converse of (3) applies with high probability \cite{Diem:CoverAttacks,Howe:IsogenyClasses}
\end{enumerate}
\end{proposition}

Additionally, the explicit isogenies attack of \cite{Smith:Genus3DiscreteLog} on genus 3 hyperelliptic curves can (and should) be avoided by ensuring $f(x)$ has exactly one irreducible factor of degree 3, 5, or 7 over $\Fq$.
In light of (1), we consider the security of a genus 3 curve $C$ over $\Fp$ comparable to that of a genus 2 curve $C'$ over $\Fpp$ when\footnote{We use ``$\lg$" to denote binary logarithms and ``$\log$" for the natural logarithm (the distinction is immaterial here and in big-O notation, but may be relevant elsewhere).}
$$\lg\#J(C/\Fp)/\lg\#J(C'/\Fpp) = 9/8$$
(note $p$ and $p'$ are different primes), and if $C$ is a genus 2 curve over $\Fptwo$ we require a ratio of 14/9.  Considering (2) and (3), for the trace zero variety $T(C/\Fpt)$, the comparable ratio is 5/4 when 3 divides $\#J(C/\Fptwo)$ and 6/5 otherwise.

See \cite{Cohen:HECHECC} and \cite{Koblitz:AlgebraicCrypto}  for further background on hyperelliptic curve cryptography.

\section{Algorithms}

To compute $\#J(C)$, we apply a probabilistic generic algorithm that computes the order (and structure) of an arbitrary finite abelian group $G$.  We give an overview of the general algorithm presented in \cite{Sutherland:Thesis}, focusing on the components most relevant to searching a family of Jacobians.

By convention, we use multiplicative notation for generic groups and let $\identity$ denote the group identity.  For $\alpha\in G$, we let $|\alpha|$ denote the \emph{order} of $\alpha$, the least positive integer $n$ for which $\alpha^n = \identity$, and call any multiple of $|\alpha|$ an \emph{exponent} of $\alpha$.  The order of $G$ is denoted $|G|$, and $\lambda(G)$ is the exponent of $G$, the least positive integer which is an exponent of every element in $G$.  Complexity metrics for generic algorithms count group operations (time) and group elements (space).

We assume the availability of a black box that uniquely identifies group elements, a requirement easily met by Jacobian arithmetic based on a Cantor-Mumford representation \cite{Cantor:Jacobian}. We also require access to randomly generated group elements, which may be obtained via the methods detailed in \cite[14.1-2]{Cohen:HECHECC}.  We presuppose a uniform distribution, however this assumption can typically be relaxed in practice.\footnote{For small groups, the simplest approach in \cite{Cohen:HECHECC} may not generate the entire Jacobian.  Decompression techniques ensure a uniform distribution (at slightly greater cost).} 

Central to our method for obtaining subexponential performance is the notion of a \emph{conditional algorithm}.  Such an algorithm is given the option to explicitly reject inputs that fail to satisfy a specified condition, but must otherwise behave correctly.%  The condition is typically one that simplifies the computation.

The condition we use here is based on the fact that for any group element $\alpha$,
\begin{equation}\label{equation:OrderReduction}
\left|\alpha^E\right| = \frac{\orda}{\gcd(\orda,E)}.
\end{equation}
This may be used to reduce the size of $|\alpha^E|$ relative to $\orda$, making it easier to compute.  This motivates the following definition.

\begin{definition}\label{definition:Easy}
Let $E=\prod q$ where $q$ ranges over maximal prime-powers bounded by $B$.  A positive integer $N$ for which $N/\gcd(N,E) \le B^2$ is called $B$-\emph{easy}, and otherwise, $B$-\emph{hard}.
\end{definition}

A $B$-easy integer is semismooth with respect to $B^2$ and $B$ (prime factors bounded by $B^2$ and all but one by $B$).  The converse holds if we apply the semismooth criteria to prime-power factors.

\begin{proposition}\label{proposition:OrderAlgorithm}
For any $B>0$ there is a probabilistic generic algorithm $\mathcal{A}$ such that the following hold for all finite abelian groups $G$:
\begin{enumerate}
\item
If $\mathcal{A}$ rejects, then $|G|$ is $B$-hard; otherwise, with high probability,\\$\mathcal{A}$ outputs $\lambda(G)$, the structure of $G$, and $|G|$.
\item 
The expected number of group operations is $O(B) + O(\lg^2|G|)$.
\end{enumerate}
\end{proposition}

We say $\mathcal{A}$ computes $|G|$ on the condition that $|G|$ is $B$-easy.
The probability and expectation in Proposition \ref{proposition:OrderAlgorithm} do not depend on $G$.  By the structure of $G$, we mean an explicit factorization of $G$ into cyclic subgroups of known order with a generator for each factor.  The structure of $G$ is computed using $\lambda(G)$, to obtain $|G|$.

When tight bounds on $|G|$ are known, as in the case of Jacobians, the fact that $\lambda(G)$ divides $|G|$ may suffice to determine $|G|$, and the structure of $G$ is not needed.  The time to compute $|G|$ is typically dominated by the time to compute $\lambda(G)$ in any event; computing the group structure at worst increases the constant factors.

If $B \ge |G|^{1/2}$ then $|G|$ is certainly $B$-easy.  By starting with a small value of $B$ and increasing it in stages, one obtains an $O(N^{1/2})$ algorithm for computing the structure of any finite abelian group.  When $|G|$ is a random integer, the median complexity is $O(N^{0.344})$.  This approach can be much faster than other generic algorithms for computing abelian group structure \cite{Buchmann:BabyGiant,Buchmann:GroupStructure,Teske:GroupStructure,Teske:PohligHellmanStructure}.

To search for Jacobians, we estimate $|G| \approx 2^n$ based on the Weil interval (\ref{equation:WeilInterval}), then pick a fixed $B = 2^{n/u}$ that minimizes $B/\sigma(u)$, where $\sigma(u)$ estimates the probability that a random integer $x$ is $x^{1/u}$-easy.  Asymptotically, we may use $\sigma(u) = G(1/u,2/u)$, where $G(s,t)$ is the semismooth probability function defined by Bach and Peralta in \cite{Bach:Semismooth} (the impact of prime-power factors is  negligible).  This yields an $L(1/2,\sqrt{2}/2)$ bound on both $1/\sigma(u)$ and $B$, leading to an $L(1/2,\sqrt{2})$ bound on the entire search, based on the heuristic assumption that $|G|$ has the distribution of a random integer of comparable size.\footnote{One uses $\sigma(u) \ge \rho(u) = u^{-u+o(1)}$, where $\rho(u)$ is the Dickman function \cite{Bach:Semismooth,Canfield:SmoothBound,Dickman:DickmanFunction}.}

It remains to prove Proposition \ref{proposition:OrderAlgorithm}.  We assume henceforth that $B \ge \lg^2|G|$.

%The simplest way for an algorithm to satisfy Proposition \ref{proposition:OrderAlgorithm} above is to reject as often as possible.  Indeed, the entire process may be viewed as an attempt to prove that the unknown quantity $|G|$ is $B$-hard, which, when unsuccessful, produces the value of $|G|$, along with the means to verify it, as a counterexample.  A high-level description of the algorithm is given below.

\begin{algorithm}[Group Order]\label{algorithm:GroupOrder}
Given a finite abelian group $G$ and a bound $B$, on the condition that $|G|$ is $B$-easy:
\vspace{3pt}
\rm
\renewcommand\labelenumi{\theenumi.}
\begin{enumerate}
\item
Compute $\lambda(G)$ on the condition that $\lambda(G)$ is $B$-easy.
\item
\vspace{6pt}
For each prime $p|\lambda(G)$, compute the structure of $H_p$, the $p$-Sylow subgroup of $G$, on the condition that $|H_p| \le B^2$. (If $|G|$ is found $B$-hard, reject).
\end{enumerate}
\it
\vspace{3pt}
Output $\lambda(G)$, the structure of $G$, and $|G|$.
\end{algorithm}

The algorithm for Step 2 is described in detail in \cite[Algorithm 9.1]{Sutherland:Thesis}.  Computation in $H_p$ is accomplished via the black box for $G$, exponentiating by $\lambda(G)/p^h$ to obtain random elements of $H_p$ (here $p^h$ is the largest power of $p$ dividing $\lambda(G)$).  Computing the structure of $H_p$ uses $O(\sqrt{|H_p|})$ group operations \cite[Proposition 9.3]{Sutherland:Thesis}.
The condition $|H_p| \le B^2$ gives a complexity of $O(B\lg|G|)$ (entirely acceptable in practice), however to prove Proposition \ref{proposition:OrderAlgorithm}, we reject whenever $|G|$ is $B$-hard, assuring an $O(B)$ bound.  As noted above, Step 2 need not be implemented to search for large Jacobians.

We now present the algorithm to compute $\lambda(G)$, using two generic subroutines that compute the order of a group element.

\begin{algorithm}[Group Exponent]\label{algorithm:Exponent}
Given a finite abelian group $G$, a bound $B$, and a constant $c$, let $E = \prod q$, where $q\le B$ ranges over prime powers, and set $N \leftarrow 1$.  On the condition that $\lambda(G)$ is $B$-easy:
\vspace{3pt}
\rm
\renewcommand\labelenumi{\theenumi.}
\begin{enumerate}
\item
For a random $\alpha\in G$, set $\alpha\leftarrow \alpha^N$, then compute $\beta \leftarrow \alpha^{E}$.
\item
\vspace{6pt}
Compute $N' = |\beta|$ on the condition that $|\beta| \le B^2$.
\item
\vspace{6pt}
Compute $N'' = |\alpha^{N'}|$ with exponent $E$.  Set $N \leftarrow NN'N''$ and $t \leftarrow 1$.
\item
\vspace{6pt}
For a random $\alpha\in G$, attempt to compute $N' = |\alpha^N|$ with exponent $E$.\\
If a reject occurs, goto Step 1.
\item
\vspace{6pt}
Set $N \leftarrow NN'$ and increment $t$.  If $t < c$ then goto Step 4.
\end{enumerate}
\it
\vspace{3pt}
Output $\lambda(G) = N$.
\end{algorithm}

The order computation in Step 2 is a bounded search for $|\beta|\le B^2$, which may be performed by a standard birthday-paradox algorithm or by Algorithm \ref{algorithm:PrimorialSteps} below.  The order computations in Steps 3 and 4 use Algorithm \ref{algorithm:LinearFastOrder}.  In Step 3, $E$ is necessarily an exponent of $\alpha^{N'}$, since $\alpha^{N'E} = \beta^{N'} = \identity$.  In Step 4, $E$ might not be an exponent of $\alpha^N$, causing Algorithm 3 to reject.  We now show this rarely happens.

If Step 2 rejects, then $|\alpha|$ is $B$-hard and so is $\lambda(G)$.  Otherwise, we claim that the integer $\lambda(G)/|\alpha|$ is $B$-powersmooth, with probability greater than $1-1/B$.  For a uniformly random $\alpha\in G$ and a prime power $p^h$ dividing $\lambda(G)$, we find that
\begin{equation}\label{equation:ExponentProbability}
\Pr\left[p^h\thickspace{\rm divides}\thickspace\frac{\lambda(G)}{|\alpha|}\right] \le \frac{1}{p^h},
\end{equation}
by considering a factorization of $G$ into cyclic subgroups of prime-power order.  If $\lambda(G)/|\alpha|$ is not $B$-powersmooth, it must be divisible by a prime power greater than $B$, and our claim is proven by (\ref{equation:ExponentProbability}).  It follows that the expected complexity is $O(B)$ group operations, assuming an $O(B)$ bound on each step.

The output value $N$ is the least common multiple of the orders of $c$ random group elements--necessarily a divisor of $\lambda(G)$.  One can apply \ref{equation:ExponentProbability} to obtain
\begin{equation}
\Pr\left[N\ne\lambda(G)\right]< 1-\frac{1}{2^{c-2}},
\end{equation}
as shown in \cite[Proposition 8.3]{Sutherland:Thesis}.

From the prime number theorem we find $\lg E \approx (B/\log B)\lg B = B/\log 2$, hence Step 1 can be accomplished with standard exponentiation techniques using $B/\log 2 + o(B)$ group operations.  Step 2 can be performed by either of the generic birthday-paradox algorithms using an expected $O(B)$ group operations.  We give a more efficient method (Algorithm \ref{algorithm:PrimorialSteps}), but it is not needed for Proposition \ref{proposition:OrderAlgorithm}.

To complete the demonstration of the proposition, we need only show how to compute $|\alpha|$ with an exponent $E$, using $O(\lg E)$ group operations.

\begin{algorithm}[Linear Order]\label{algorithm:LinearFastOrder}
Given an integer $E=q_1q_2\cdots q_w$ factored into prime powers and $\alpha\in G$, let $\alpha_0=\alpha$.  On the condition that $\alpha^E = \identity$:
\vspace{3pt}
\rm
\begin{enumerate}
\item
For $i=1$ to $w$, compute $\alpha_{i} \leftarrow \alpha_{i-1}^{q_i}$ until $\alpha_{i} = \identity$.\\
If this fails to occur then reject, otherwise set $N\leftarrow |\alpha_{i-1}|$.
\vspace{6pt}
\item
Do a binary search for the least $j\in[0,i]$ for which $\alpha_j^N = \identity$.\\
If $j > 0$ then set $i \leftarrow j-1$, $N \leftarrow |\alpha_{i}^N|N$, and repeat Step 2.
\end{enumerate}
\vspace{3pt}
\it
Output $\orda = N$.
\end{algorithm}
The correctness of the algorithm follows from the invariant $N = |\alpha_i|$.
Each of the assignments to $N$ involves a prime-power order computation accomplished via repeated exponentiation.  The cost of Step 1 is at most $m + o(m)$ group operations, where $m = \lg E$.  The cost of Step 2 may be bounded by $O\bigl((n^2/\lg n)\lg m\bigr)$, where $n = \lg|\alpha|\le \lg|G|$.  For $m \ge n^2$, the total complexity is $O(m) = O(\lg E)$.  For $m \gg n^2$ it is at most $m + o(m)$, possibly much less.\footnote{In Step 4 of Algorithm 4, $\orda$ quickly becomes very smooth.}

This completes our demonstration of Proposition \ref{proposition:OrderAlgorithm}.

\begin{proposition}\label{proposition:OrderAlgorithmSpace}
Algorithm $\mathcal{A}$ of Proposition \ref{proposition:OrderAlgorithm} can be implemented using storage for $O(\lg^2|G|)$ group elements.
\end{proposition}
\begin{proof}[Proof sketch]
Pollard's rho method can be used in Step 2 of Algorithm \ref{algorithm:GroupOrder} and in Step 2 of Algorithm \ref{algorithm:Exponent}, using storage for $O(\log|G|)$ group elements.  This also suffices for all exponentiations ($E$ need not be explicitly computed).  As written, Algorithm \ref{algorithm:LinearFastOrder} stores $O(w) = O(B/\log B)$ group elements.  When $B/\log B > \lg^2|G|$, the space can be made $O(\lg^2|G|)$ by saving only this many values of $\alpha_i$, recomputing as required, with a negligible impact on time (see \cite[Proposition 7.1]{Sutherland:Thesis} for details).
\end{proof}

When searching for $B$-easy groups, it is typically not necessary to constrain space to this extent.  The chosen bound $B$ is $L(1/2,\sqrt{2}/2)$ in terms of $N=|G|$ (unconditionally), and the space will be $O(B)$ (or better) even if a rho search is not used.  The algorithm described below is significantly faster than a rho search, and in practice the space requirements are moderate (see Section \ref{section:Examples}).

\subsection*{A Parallel Primorial-Steps Algorithm}
The primorial-steps algorithm \cite{Sutherland:Thesis} computes the order of an element in a generic group.  It is asymptotically faster than the standard $\Theta(\sqrt{N})$ birthday-paradox algorithms, with an improvement of $\Theta(\sqrt{\log\log N})$.  In practice the gain is a factor of two or three in both time and space over a standard baby-steps giant-steps implementation.  We present a parallel version of the algorithm, designed for black boxes that can perform parallel group operations more quickly (typically by combining inversions in an underlying field).

A parallel algorithm is well suited to a multi-processor environment, but our primary motivation here is to speed up the group operation in a single thread of execution. The resulting algorithm remains generic.  The performance improvement depends on the black box, but can be substantial (see Table \ref{table:BlackBoxes}).

As above, our approach is based on (\ref{equation:OrderReduction}).  Given $\alpha\in G$, we use exponentiation by a suitable $E$ to remove small primes from $|\alpha^E|$, enabling an optimized baby-steps giant-steps search.  When $B = 1000$, we use $E = 2^{19}\cdot 3^{12}\cdot 5^8\cdot 7^7$, and compute $\beta = \alpha^E$. Assuming that $\orda \le B^2$, the integer $N=|\beta|$ is then relatively prime to the \emph{primorial} $P_4 = 2\cdot 3\cdot 5\cdot 7 = 210$.  In general, we choose a primorial $P = P_w$ and a positive integer $m$, so that
\begin{equation}
m^2P\varphi(P) \ge B^2,\notag
\end{equation}
where $\varphi(P)$ is Euler's function.  We compute $m\varphi(P)$ baby steps $\beta^b$ for each $b\in[1,mP]$ relatively prime to $P_w$, and a similar number of giant steps $\beta^{mPa}$ for $a$ from 1 to $\varphi(P)$.
Since any integer $N \le B^2$ may be written in the form
\begin{equation}
N = mPa - b,\notag
\end{equation}
one of our baby steps must match one of our giant steps.
In our example, we let $m = 10$ and use 480 baby steps followed by 477 giant steps.  We require some additional group operations to compute $\beta = \alpha^E$ and the values $\beta^2$, $\beta^4$, \ldots, $\beta^{10}$ needed to span the gaps between integers relatively prime to $P_4$.  The total is about half the 2000 group operations used in a standard baby-steps giant-steps search, and the space is similarly reduced.  In general the improvement is a factor of $\sqrt{P/\varphi(P)}$, and for a suitable $P$, one obtains an asymptotic complexity of $O(B/\sqrt{\log\log B})$.\footnote{This generalizes to an unbounded search for $N=\orda$ with complexity $O(\sqrt{N/\log\log N})$ \cite{Sutherland:Thesis}.}

Applying optimized baby-step giant-step methods given constraints on $|\beta|$ is not new; this technique is often used in conjunction with $\ell$-adic methods, as in \cite{Matsuo:BabyStepPointCounting}.  The novelty here is that the constraints are obtained generically.

In the example above we could have used $P_5$ rather than $P_4$, and set $m$ to 1.  We intentionally use a slightly suboptimal value of $P$ and a larger value of $m$ to facilitate a parallel implementation.  Rather than a single sequence of baby steps, we use $m$ sequences, each spanning a range of $P$ powers of $\beta$ using $\varphi(P)$ group operations.  We similarly use $m$ suitably spaced sequences of $\varphi(P)$ giant steps.

To incorporate the usual optimization for fast inverses, we double the spacing between giant steps.  We now choose $m$ and $P = P_w$ satisfying
\begin{equation}
2m^2P\varphi(P) \ge B^2,\notag
\end{equation}
and assume these values are precomputed, along with the value $E = \prod_{i=1}^w p_i^{h_i}$, where $p_i$ is the $i$th prime and $p_i^{h_i}\le B^2 < p_i^{h_i+1}$.

We also precompute a \emph{wheel} for the primorial $P$, a sequence $r(n)$ with the property that $1+\sum_{j=1}^nr(j)$ gives the $(n+1)$st positive integer relatively prime to $P$, for $n$ from 1 to $\varphi(P)$.  The wheel for $P_3 = 30$ is the sequence (6,4,2,4,2,4,6,2).  The value $r_{\rm max}$ denotes the largest element of the sequence $r(n)$.  To obtain greater flexibility in the choice of $m$, one can use $P=tP_w$ and ``roll" the wheel $t$ times in Step 3 below.  For simplicity we assume $P = P_w$.

\begin{algorithm}[Primorial-Steps]\label{algorithm:PrimorialSteps}
Given $\alpha \in G$ and a bound $B$, let $m$, $P$, $E$, and $r(n)$ be as above, and let $\mathcal{B}$ and $\mathcal{G}$ be empty sets. On the condition $|\alpha| \le B^2$:
\vspace{3pt}
\rm
\renewcommand\labelenumi{\theenumi.}
\begin{enumerate}
\item
Compute $\beta_0\leftarrow\alpha^E$ and $\delta_i = \beta_0^{i}$, for even $i$ from 2 to $r_{\rm max}$.\\
Compute $\beta_i \leftarrow \beta_0^{P}\beta_{i-1}$, for $i$ from 1 to $m-1$.  Let $\vec{\beta} = (\beta_0,\ldots,\beta_{m-1})$.
\vspace{6pt}
\item
Compute $\gamma_0 \leftarrow \beta_0^{mP}$ and $\delta_0 = \gamma_0^2$.\\
Compute $\gamma_i = \gamma_0^{2\varphi(P)}\gamma_{i-1}$, for $i$ from 1 to $m-1$.  Let $\vec{\gamma} = (\gamma_0,\ldots,\gamma_{m-1})$.
\vspace{6pt}
\item
Set $k \leftarrow 1$ and for $i$ from 0 to $m-1$, set $\mathcal{B}\leftarrow\mathcal{B}\cup(\beta_i, i, k)$.\\
For $j$ from 1 to $\varphi(P)-1$: compute $\vec{\beta}\leftarrow\vec{\beta}\delta_{r(j)}$, set $k \leftarrow k + r(j)$,\\ \indent and for $i$ from 0 to $m-1$, set $\mathcal{B}\leftarrow\mathcal{B}\cup(\beta_i, i, k)$.
\vspace{6pt}
\item
If $(\identity, i, k)\in \mathcal{B}$, minimize $N = Pi+k$ over such tuples and goto Step 6.
\vspace{6pt}
\item
Set $k \leftarrow 0$, then for $i$ from 0 to $m-1$ set $\mathcal{G}\leftarrow\mathcal{G}\cup(\gamma_i, i, k)$.\\
For $k$ from 1 to $\varphi(P)-1$: compute $\vec{\gamma}\leftarrow\vec{\gamma}\delta_0$,\\
\indent and for $i$ from 0 to $m-1$, set $\mathcal{G}\leftarrow\mathcal{G}\cup(\gamma_i, i, k)$
\vspace{6pt}
\item
Find the least $N = a \pm b$ corresponding to $\gamma\beta^{\pm 1} = \identity$,\\
\indent where $\thickspace a = 2mP(\varphi(P)i_1 + k_1)$ and $b = Pi_2+k_2$, \\
\indent for some $(\gamma, i_1, k_1)\in \mathcal{G}$ and $(\beta, i_2, k_2)\in \mathcal{B}$.\\
If no such $N$ exists, reject.
\vspace{6pt}
\item
Compute $N' = |\alpha^N|$ with exponent $E$ and set $N\leftarrow NN'$.
\end{enumerate}
\vspace{3pt}
\it
Output $N = \orda$.
\end{algorithm}

In a parallel search, the first match found doesn't necessarily give the order of $\beta_0$ (we could find a multiple of $|\beta_0|$), hence the minimization of $N$.  In a standard implementation we don't explicitly construct the set $\mathcal{G}$, rather we check for a match in $\mathcal{B}$ as each giant step is computed.  This allows early termination when successful, provided we handle the case that $N$ is a multiple of $|\beta_0|$.  We can determine $|\beta_0|$ using the (factored) exponent $N$ before proceeding to Step 7.  

In the present application there is good reason to explicitly compute $\mathcal{G}$ and perform the matching process in Step 6, as shown.  We expect the search to fail in most cases, so this doesn't materially impact the time.  It doubles the space required, but space is not a limiting factor and there are other ways to reduce it.  The tuples in $\mathcal{B}$ and $\mathcal{G}$ needn't store entire group elements (a small hash value suffices) and $B$ may be chosen judiciously (see Table \ref{table:Semismooth}).\footnote{In fact, explicit computation of $\mathcal{B}$ and $\mathcal{G}$ enables  searches that are effectively unlimited by RAM--they can be efficiently migrated to disk and matched via a merge or radix sort.}

The advantage of computing both $\mathcal{B}$ and $\mathcal{G}$ explicitly is that it allows matching to be performed more efficiently (it also enables greater parallelism).  When the group operation is extremely fast, the implementation of the lookup table used in a baby-steps giant-steps search can have a significant impact on the running time.  On our test platform, the fastest black boxes achieve execution times close to the latency of general memory access.  If we defer matching until the end it can be done more quickly, with better locality of reference, as described in Section \ref{section:Implementation}.

\subsection{Recovering the zeta function}\label{subsection:ZetaFunctionRecovery}
Having computed $\#J(C)$, we need to determine the zeta function of $C$.  This problem (and many others) is discussed in \cite{Elkies:CurvesOverFiniteFields}.  We provide explicit details here for the genus 2 and 3 cases and analyze the cost of determining the zeta function once $\#J(C)$ is known.
\begin{lemma}
Let $P(z)$ denote the $L$-polynomial of a non-singular, irreducible, projective curve $C$ of genus $g\le 3$ defined over $\Fq$.  For sufficiently large $q$, the values $P(1)$ and $P(-1)$ uniquely determine the coefficients of $P(z) = \sum_{i=0}^{2g}a_iz^i$.
\end{lemma}
\begin{proof}
Recall that $a_0 = 1$ and $a_{2g-i} = q^{g-i}a_i$ for $0\le i < g$ (Theorem \ref{theorem:Weil}).
If $g\le 2$ then $a_1 = \left[P(1)-P(-1)\right]/[2(q+1)]$, and for $2\le g\le 3$, we find that $a_2 = \left[P(1)+P(-1)-2(p^2+1)\right]/2$.  This proves the lemma for $g\le 2$.  For $g = 3$,
$$a_1 = \frac{P(1) - P(-1) + 2a_3}{2(q^2+1)} = A_1 + \delta_1,$$
where $A_1$ is fixed and $|\delta_1| < \binom{6}{3}q^{-1/2}$, by the bounds in (\ref{equation:abounds}).  If $q \ge 40^2$ then $|\delta_1| < 1/2$ and the integer $a_1$ is determined, fixing $a_3$ as well.
\end{proof}
To determine $P(-1)$, we compute $\#J_{2/1}(C)=P(1)P(-1)/P(1)$ (Lemma \ref{lemma:ZetaCount}).  For hyperelliptic curves, we may equivalently compute $\#J(\Ctwist)$ (Lemma \ref{lemma:Twist}).  In the non-hyperelliptic case we can compute in $J_{2/1}(C)$ via group operations in $J_2(C)$, using exponentiation by $P(1)$ to obtain elements of $J_{2/1}(C)$.\footnote{If $d = \gcd(P(1),P(-1)) > 1$, we lose $d$-torsion elements of $J_{2/1}(C)$ when we do this, but the resulting subgroup will usually be large enough for our purposes.}

For simplicity, we assume $C$ is hyperelliptic.  In genus 2 we have
\begin{equation}\notag
a_1 = \frac{P(1) - (q^2+1) - a_2}{q+1} = \frac{P(1)-(q^2+1)}{q+1} - \frac{a_2}{q+1} = A_1 + \delta_1,
\end{equation}
where the bounds in (\ref{equation:abounds}) imply $|\delta_1| < \binom{4}{2}$.  There are at most eleven possible values for the integer $a_1$. The corresponding values $P(-1) = P(1) - 2(q+1)a_1$ form an arithmetic sequence with difference $2(q+1)$.  To distinguish the correct $P(-1)$, we generate a random $\alpha\in J(\Ctwist)$ and step through the sequence, computing  powers of $\alpha$.  If only one exponent of $\alpha$ is found then we have determined $P(-1)$ (usually the case).  Otherwise we compute $\orda$ using a suitable exponent (e.g., the gcd of all exponents found).  By repeating this process we determine $\lambda(J(\Ctwist))$ and, if necessary, apply Step 2 of Algorithm \ref{algorithm:GroupOrder} to compute $\#J(\Ctwist)=P(-1)$.

We use a similar approach in genus 3, and assume $q > 1640$.  We then have
\begin{equation}\notag
a_1 =  \frac{P(1) - (q^3+1) - (q+1)a_2 - a_3}{q^2+1} = A_1 + \delta_1,
\end{equation}
with $|\delta_1| < \binom{6}{2} + 1/2$, giving 31 possibilities for $a_1$.  Given $P(1)$ and $a_1$, we find
$$a_2 = \frac{P(1) - (q^3+1) - (q^2+1)a_1 - a_3}{(q+1)} = A_2 + \delta_2,$$
with $|\delta_2| < \binom{6}{3}q^{1/2}$.  The corresponding $P(-1)$ values are given by
$$P(-1) = 2(q^3+1) - P(1) + 2(q+1)a_2,$$
and we now have 31 arithmetic sequences for $P(-1)$, each with a difference of $2(q+1)$ and length less than $40q^{1/2}$.  For a random $\alpha\in J(\Ctwist)$ we can search for exponents of $\alpha$ among all these sequences simultaneously using a baby-steps giant-steps search (but not a primorial-steps search).  We compute roughly $\sqrt{1240q^{1/2}}$ consecutive powers of $\beta = \alpha^{2(q+1)}$ (baby steps), followed by a similar number of giant steps suitably spaced among the 31 sequences.  This can be implemented using parallel group operations, as in Algorithm \ref{algorithm:PrimorialSteps}.

The total number of group operations is $O(q^{1/4})$ in genus 3.  In terms of the group size $N \approx q^g$, this is $O(N^{1/12})$, with a leading constant factor of about 70.  Given the heuristically subexponential time required to find $\#J(C)$, the $O(N^{1/12})$ term is asymptotically dominant, but for practical values of $N$ this is not the case.  Even for $N$ as large as $2^{256}$, given $\#J(C)$, the time to recover the zeta function is entirely tractable in genus 3 (perhaps a few minutes).

With modification, this approach can be applied in higher genera, but the running time becomes more significant.  In genus 4 the complexity is $O(N^{3/16})$ group operations, and in general the complexity is $\Omega(N^{(g-1)(g-2)/(8g)})$.  Extensions to Kedlaya's algorithm \cite{Bostan:LinearRecurrences,Harvey:LargeCharacteristicKedlaya} compute the zeta function of a curve in $\tilde{O}(p^{1/2})$ time, giving an $\tilde{O}(N^{1/(2g)})$ algorithm.  This is slightly faster than our method in genus 4, and the advantage grows quickly in higher genera.

\section{Implementation}\label{section:Implementation}

We mention a few implementation details that may be relevant to those wishing to replicate our results.
Our implementation platform was a 2.5GHz AMD Athlon 64 processor (dual core) with 2GB of memory, running a 64-bit Linux operating system.  We ran eight of these systems in parallel in the larger tests.  The algorithms were implemented using the GNU C compiler \cite{GNU} and the GMP multi-precision arithmetic library \cite{GMP}.

\subsection*{Black Boxes}

The parallel group operation enabled by Algorithm \ref{algorithm:PrimorialSteps} is most advantageous to a black box based on an affine representation of the Jacobian.  We used modified versions of Algorithms 14.19-21 (genus 2) and Algorithms 14.52-53 (genus 3) in \cite{Cohen:HECHECC}.  The black box executes several group operations up to the point where a field inversion is required, performs a single combined field inversion using Montgomery's trick \cite{Cohen:CANT}, then completes the group operations.  With this approach the amortized cost of a field inversion is 3 field multiplications (3M), and the effective cost of a group operation is then 28M in genus 2 and 74M in genus 3.
%This is substantially faster than alternative representations.

%We implemented the prime field arithmetic using either a one or two word 64-bit representation.  All the code was 
The prime field arithmetic was implemented in C except for two in-line assembly directives, one to compute the 128-bit product of two 64-bit values and one to perform a 128-bit addition.  The field multiplications for the Mersenne primes $2^{61}-1$ and $2^{89}-1$ were specifically optimized, but otherwise we used a Montgomery representation \cite{Montgomery:ModularArithmetic}, and Montgomery inversion was used in all cases (see \cite{Cohen:HECHECC} and \cite{Menezes:Handbook} for algorithms).  Performance metrics appear in Table \ref{table:BlackBoxes}.

\subsection*{Parallel exponentiation}

Both asymptotically and in practice, the exponentiation performed in Step 1 of Algorithm \ref{algorithm:Exponent} dominates the total running time.  To obtain the performance improvements offered by parallel group operations, we exponentiate in parallel for several curves defined over the same field.\footnote{For families of curves where the field varies, one might include 100 curves per field.}  The exponent $E$ does not change once $B$ is chosen, so it should be precomputed and put into a convenient form (a $2^k$-ary sliding representation was used in our tests).  The exponentiation takes more than twice as long as the search step, so it is convenient to have two threads performing exponentiations on a dual processor, feeding their results to a single search thread running in parallel.

\subsection*{Choosing the bound $\boldsymbol{B=2^{n/u}}$}

Given the comments above, we might search to a bound greater than $B^2$, balancing the time between exponentiation and searching.  The behavior of the semismooth probability function $G(1/u,1/v)$ argues against this, as small changes in $v$ have little impact.  In fact, the optimal choice of $v$ is slightly above $2/u$, implying a bound less than $B^2$, but the difference is negligible.

The quantity $2^{n/u}/\sigma(u)$ is insensitive to small changes in $u$ close to the optimal value.  We can choose a slightly smaller $u$ without materially impacting the running time, obtaining a substantially smaller $B=2^{n/u}$, which saves space.  As seen in Table \ref{table:Semismooth}, one can reduce space by more than a factor of two while increasing the time by only $5\%$.
%The fact that we generally find our actually success rate to be better than predicted by $G(1/u,2/u)$ (Table \ref{table:Actuals}) also motivates a slightly larger value of $u$.

\subsection*{Eliminating 2-torsion}

One can efficiently filter a family of hyperelliptic curves of the form $y^2 = f(x)$ to remove curves whose Jacobian has even order by testing whether $f(x)$ is an irreducible polynomial in $\Fp[x]$.  This is well worth doing if one is interested in the group $J_{2/1}(C) = J(\Ctwist)$, since $\#J(C)$ and $J\#(\Ctwist)$ must have the same parity, but otherwise the situation is less clear.  As shown in Table \ref{table:Distribution}, while the probability of finding groups with near-prime order generally increases when $\#J(C)$ is odd, the probability that $\#J(C)$ is $B$-easy goes down, more than offsetting the increase in many cases.  Note that it is possible for any of the groups $J_{4/2}(C)$, $J_{3/1}(C)$ and $J_{3/1}(\Ctwist)$ to have prime order even when $\#J(C)$ is even.

\subsection*{Efficient matching}

In the description of Algorithm \ref{algorithm:PrimorialSteps}, the sets $\mathcal{B}$ and $\mathcal{G}$ contain tuples $(\alpha, i, k)$, where $\alpha\in G$ is a baby-step or a giant-step.  It is not necessary to store $\alpha$.  It can be recovered using $i$ and $k$ by exponentiating $\beta_0$ or $\gamma_0$, which are known to the algorithm. Asymptotically, a uniform hash value of $\lg B (\lg\lg B)^{1+\epsilon}$ bits suffices to keep the cost of matching negligible.  In practice, $\lg B$ is less than 30, the values $i$ and $k$ require a total of $\lg B$ bits, and we use a $(64-\lg B)$-bit hash value to make a 64-bit value for each tuple.  This is about one third the size of a compressed group element.

When using a fast inverse optimization, it is helpful if an element and its inverse hash to the same value.  This allows detection of both the cases $\gamma\beta = \identity$ and $\gamma\beta^{-1} = \identity$, without requiring extra table entries.  To localize memory access, a merge sort or (better) a radix sort may be used to find matches \cite[5.2.4-5]{Knuth:AOPIII}.  We used a partial radix sort with a radix of $2^9$ or $2^{10}$.
%This approach can be applied to any baby-steps giant-steps implementation.  When using a fast inverse optimization, as in Algorithm \ref{algorithm:PrimorialSteps}, it is helpful if an element and its inverse hash to the same value (easily accomplished with Jacobians).  To facilitate matching with localized memory access, we do a partial radix sort.  We initially partition the sets into $2^m$ linear lists based on the low order $m$ bits of the hash value at the time of insertion, where $m$ is small enough to keep the end of each list in cache (we used $m=9$).  Once the sets are complete, we further partition corresponding lists until we reach a conveniently small size to find matching hash values efficiently.  For each matching pair of hash values we then reconstruct the corresponding group elements to test whether $\gamma\beta^{\pm 1} = \identity$.

\section{Examples}\label{section:Examples}

%The process of searching for cryptographically suitable Jacobians is perhaps best illustrated by example, sowe present these in some detail before giving the algorithms.  The few mathematical facts we need are well known and can be found in Section \ref{section:Facts}.
For ease of illustration, we use parameterized families of curves with small coefficients.  This choice is arbitrary, as is the choice of finite field.  In practice, one might choose a family of curves whose coefficients admit a particularly efficient implementation, as suggested by Bernstein \cite{Bernstein:EllipticVsHyperelliptic}.
\subsection*{Genus 3 examples}
We use the family of hyperelliptic curves defined by
$$y^2 = x^7 + 3x^5 + x^4 + 4x^3 + x^2 + 5x + t$$
over the prime field $\Fp$, with $p = 2^{50}-27$.  As we are interested in groups $J(C)$ with near-prime order, we don't try to compute $\#J(C)$ directly.  Instead, we attempt to compute $\#J(\Ctwist)$ for each curve in our family, using Algorithm 1 (with random coefficients, twisting is unnecessary).  When successful, we apply the method described in Section \ref{subsection:ZetaFunctionRecovery} to recover the $L$-polynomial of $\Ctwist$ using $\tilde{P}(1) = \#J(\Ctwist)$. We then compute $\#J(C) = \tilde{P}(-1) = P(1)$.

Given the bound $B$, Algorithm 1 succeeds if $\#J(\Ctwist)$ is $B$-easy (Definition \ref{definition:Easy}).  To choose $B$, we make the heuristic assumption that $\#J(\Ctwist)$ is a random $n$-bit integer, where $n = \lg\#J(\Ctwist) \approx 150$. We pick $u$ to minimize $2^{n/u}/\sigma(u)$ and set $B = 2^{n/u}$.  For $n = 150$, we choose $u = 6.25$ and find $\sigma(u) \approx 1/1765.$ (Table \ref{table:Semismooth}).

Algorithm 1 uses $B/\log 2 + o(B)$ group operations, about 36 million in this case (Table \ref{table:Actuals}).  The genus 3 black box performs roughly 1.8 million group operations per second (Table \ref{table:BlackBoxes}), and the CPU time per curve is about 20 seconds on a 2.5 GHz AMD Athlon-64.  This chip has two processors, so on a single PC we test a curve every 10 seconds.  We achieve our first success when $t = 648$ and find
$$\#J(\Ctwist) = 2^3\cdot 5^2 \cdot 233 \cdot 937\cdot 8053\cdot 18719\cdot 44171\cdot 1180799\cdot 13517389\cdot 307558308259.$$
Given $\#J(\Ctwist)$, it takes only 300,000 group operations (0.2 seconds) to determine the zeta function of $C$, whose $L$-polynomial $P(z)$ has coefficients
$$a_1 = 39141148,\thickspace a_2 = 1354965780525799,\thickspace a_3 = 18939879984661962930696.$$
We then compute
$$\#J(C) = P(1) = 2^3\cdot 3\cdot 1083611 \cdot 54880077749424473770842486727458448993.$$
This value isn't quite near prime, but it could have been; see Section \ref{table:Genus3Curves} for a prime example.  The average time required to successfully compute $\#J(C)$ for some curve in our family is about four hours on a single PC, somewhat better than $\sigma(u)$ would suggest (we generally find our heuristic assumption pessimistic).  The memory requirements are modest, about 200MB in this case, and in our largest tests, about 1GB.  Memory usage can be substantially reduced with a small ($<5\%$) impact on performance (see Table \ref{table:Semismooth}).\footnote{By Proposition \ref{proposition:OrderAlgorithmSpace}, the space can be made $O(\lg^2\#J(C))$, but with a larger impact on time.}    We optimized for time.

We tested similar families of curves with $p = 3\cdot 10^{16}+29$ (164-bit group size) and $p = 2^{61}-1$ (183-bit group size).  Sample results are given in Table \ref{table:Genus3Curves}.  In the first case it took about a day per success on a single PC, and in the second, slightly over four days.  We used eight PCs, succeeding roughly twice a day in the larger test.  Distributed computation not only increases the throughput, it gives a linear speedup in the time to achieve the first success for up to $O(1/\sigma(u))$ processors.  After an initial partitioning of the family of curves, no communication is required, making a distributed implementation straightforward.

It is entirely feasible to find cryptographically suitable genus 3 Jacobians using this approach.  However, the genus 2 case is more attractive, as we can find groups offering better security in much less time.

\subsection*{Genus 2 curves (1)}
We first consider the family of curves defined by
$$y^2 = x^5 + 2x^3 + 7x^2 + x + t$$
over the prime field $\Fp$ with $p = 2^{61}-1$.  We don't twist our initial family of curves, as we are not interested in the group $J_{2/1}(C) = J(\Ctwist)$, but rather the other three groups listed in (\ref{equation:candidates}).  In this case $n = \lg\#J(C) \approx 122$ and we let $u = 5.8$, obtaining $B = 2^{n/u} \approx 2^{21}$ and $\sigma(u) \approx 1/549.$

Algorithm 1 now uses about 4.6 million group operations per curve and the black box performs over 4 million group operations per second.  The CPU time per curve is under 1.1 seconds, so we test about two curves per second and expect to successfully compute $\#J(C)$ roughly every five minutes on a single PC.

For each value of $\#J(C)$, we recover $P(z)$ and, applying Lemma \ref{lemma:ZetaCount}, compute
\begin{align}\notag
\#J_{3/1}(C) &= P(\omega)P(\omega^2),\\\notag
\#J_{3/1}(\Ctwist) &= P(-\omega)P(-\omega^2),\\\notag
\#J_{4/2}(C) &= P(i)P(-i),\notag
\end{align}
where $\omega = e^{2\pi i/3}$.  If any of these are near prime, we may have found a cryptographically suitable group.\footnote{In this example $\#J_{4/2}(C)$ is not quite large enough.  We include it for the sake of illustration.}  The first case where we succeed in computing $\#J(C)$ occurs when $t = 816$, and we find that the $P(z)$ has coefficients
$a_1 = 618350030$ and $a_2 = 415833882783789026$.  The most interesting value is
$$\#J_{3/1}(C) = P(\omega)P(-\omega) =  5^2\cdot 547\cdot P_{231},$$
where $P_{231}$ is a 231-bit prime.
%P_231 = 2067243365913945149159824423420249130586442464271255167705725657330023,$$
In this case $J_{3/1}(C)$ is equal to the trace zero variety $T_3(C)$.  Even after taking Proposition \ref{proposition:Security} into account, which reduces the effective security of Jacobians over extension fields, this is well into cryptographic range.

The next interesting case occurs on our fifth success, when $t = 3909$.  We find
$$\#J_{4/2}(C) = P(i)P(-i) = 41^2\cdot P_{234}.$$
%16817104721269571216430467060464252947111108155674056433150220327309509.$$
This is the order of $J_2(\Ctwist_2)$ where $\Ctwist_2$ is the quadratic twist of $C$ in $\Fptwo$.  This curve may be written as
$$y^2 = x^5+2\alpha^2 x^3+7\alpha^3 x^2 + \alpha^4 x + 3909\alpha^5,$$
where $\alpha$ is any non-residue in $\Fptwo$.
%We have a Jacobian over $\Fptwo$ with 234-bit near-prime order, however the effective security is not quite in cryptographic range.
On our eleventh successful computation, when $t = 6005$, we find
$$\#J_{3/1}(\Ctwist) = P(-\omega)P(-\omega^2) = 4\cdot P_{242},$$
giving a group with security comparable to a 194-bit genus 2 Jacobian of prime order over a prime field.  

The total time to reach this point is about 50 minutes, a typical scenario.  If we are willing to wait a bit longer, we can find many groups with prime order, including cases where both $\#J_{3/1}(C)$ and $\#J_{3/1}(\Ctwist)$ are prime (Table \ref{table:Genus2Curves}).  On a 64-bit platform, computation in either $J_2(C)$ or a trace zero variety $T_3(C)$ with $p=2^{61}-1$ will likely be faster than in a Jacobian over a larger prime field with comparable security.

%We should note that, on 32-bit and, especially, 64-bit hardware platforms, the performance of the group operation in trace zero varieties over $\Fpt$ or Jacobians over $\Fptwo$ of the sizes considered here is likely to be superior to performance in a Jacobian over a prime field, even at the same security level.

\subsection*{Genus 2 curves (2)}
We tested genus 2 curves over three larger fields.  In the first test we used the family of curves defined by $$y^2 = x^5 + x + t$$ over the prime field $\Fp$, with $p = 2^{84}-35$.  For $u=6.5$ it takes about 34 hours per success on a single PC with $\#J(C)\approx 2^{168}$.  %This compares to one week per curve reported by Gaudry and Schost \cite{Gaudry:Genus2PointCounting} over a comparable size field.
We computed the order of 43 Jacobians, finding 11 groups with near-prime order, including:
\begin{enumerate}
\vspace{3pt}
\item
The Jacobian of the curve $y^2 = x^5+x+127861$ over $\Fp$ has near-prime order, with a 160-bit prime factor and a cofactor of 288.
\item
\vspace{3pt}
The trace zero variety of the curve $y^2 = x^5+x+89993$ over $\Fpt$ has 336-bit prime order.
%\item
%\vspace{3pt}
%The trace zero variety of the curve $y^2 = x^5+x+87058$ over $\Fpt$ has near-prime order (cofactor 31) and the trace zero variety of the quadratic twist has prime order.
\vspace{3pt}
\end{enumerate}

Our second test used the same family of curves with $p = 2^{89}-1$.  With $u=6.7$ it takes about 4 days per success per PC with $\#J(C)\approx 2^{178}$.  We computed the order of 31 Jacobians, again finding 11 groups with near-prime order, including:
\begin{enumerate}
\vspace{3pt}
\item
The Jacobian of the curve $y^2 = x^5+x+202214$ over $\Fp$ has near-prime order, with a 171-bit prime factor and a cofactor of 180.
\item
\vspace{3pt}
The Jacobian of the curve $y^2 = x^5 +\alpha^4x +207686\alpha^5$ over $\Fptwo$ has near-prime order, with a 349-bit prime factor and a cofactor of 169.
\item
\vspace{3pt}
The trace zero variety of the curve $y^2 = x^5 + 81x + 15466464$ over $\Fpt$ has near-prime order, with a 354-bit prime and a cofactor of 7.  %This group is comparable to a 295-bit genus 2 Jacobian of prime order over a prime field.
\end{enumerate}
%The security offered by the last group is comparable to a 295-bit genus 2 Jacobian of prime order over a prime field.

Our largest test used $p=2^{93}-25$.  We computed the order of a 186-bit Jacobian, finding a 372-bit trace zero variety of near-prime order for the curve
$$y^2 = x^5+2x^3+3x^2+5x+1050.$$
See Table \ref{table:Genus2Curves} for the zeta functions of all the curves mentioned above.

\section{Conclusion}\label{section:Conclusion}

For general families of genus 2 and genus 3 curves, efficiently finding cryptographically suitable Jacobians over prime fields remains a challenge.  Our method substantially increases the size of Jacobians whose order can be effectively computed, and is feasible at the low end of the cryptographic range.  In a distributed implementation, 200-bit group sizes are within reach.  As our algorithm does not use $\ell$-adic or $p$-adic methods, a combined approach may offer further improvement.  %The modular equations defined in \cite{Gaudry:ModularEquations}, for example, might be applied to obtain divisibility constraints for small primes, allowing one to increase both the probability that $\#J(\Ctwist)$ is $B$-easy and the probability that $\#J(C)$ is near-prime.

Given a family of genus 2 curves defined over a prime field, we can find cryptographically suitable groups over low degree extension fields efficiently on a single PC.  Groups offering security comparable to 200-bit genus 2 Jacobians over prime fields are easily obtained (about an hour), and the time required to achieve 250-bit security levels is not unreasonable (a day or two).  Trace zero varieties, in particular, appear to offer an attractive combination of performance and security.

%A key element of the results presented here is the use of the zeta function to connect the orders of two different groups.  This allows a hard problem to be transformed into one that may be easier, giving subexponential performance in a suitable distribution of groups.  This approach may be applicable in other instances where such a relationship can be established.

\section{Acknowledgments}
The author would like to thank Kiran Kedlaya for suggesting the problem of point-counting on hyperelliptic curves and helpful feedback on early drafts of this paper.   Thanks are also due to Rene Peralta for providing a program to compute the semismooth probability function, $G(\alpha,\beta)$.
%The source code can be found at \url{http://www.cs.yale.edu/homes/peralta/papers/semismooth.c}.

\begin{table}
\begin{center}
\begin{tabular}{@{}llrrrrr@{}}
\toprule
Black Box & $|\Fp|$ & $\lg|G|$ &\hspace{24pt} $\times 1$ & $\times 100$ &\hspace{6pt} $E$ &\hspace{6pt} $S$ \\
\midrule
Genus 2	& $2^{50}-27$ & 100 & 1.49 & 4.26 & 4.07 & 3.42\\
			& $2^{61}-1$  & 122 & 1.35 & 4.81 & 4.76 & 3.82\\
%			& $2^{75}-97$ & 150 & 0.71 & 1.80 & 1.72 & 1.60\\
%			& $2^{82}-57$	& 164 & 0.67 & 1.83 & - & -\\
			& $2^{84}-35$	& 168 & 0.66 & 1.85 & 1.78 & 1.68\\
		 	& $2^{89}-1$	& 178	& 0.67 & 2.11 & 2.01 & 1.87\\
		 	& $2^{94}-3$	& 188 & 0.62 & 1.83 & 1.76 & 1.63\\
\midrule
Genus 3	& $13^{11}+34$ & 122 & 1.12 & 1.84 & 1.86 & 1.62\\
			& $2^{50}-27$ & 150 &  1.02 & 1.83 & 1.85 & 1.67\\
			& $3\cdot 10^{16}+29$ & 164 & 0.98 & 1.83 & 1.85 & 1.67\\
			& $2^{61}-1$	& 183	& 0.92 & 1.83 & 1.85 & 1.68\\
\bottomrule
\vspace{6pt}
\end{tabular}
\caption{Black Box Performance for Prime Fields}\label{table:BlackBoxes}
\end{center}
\begin{minipage}{1.0\linewidth}
\small
The last four columns list average performance in millions of group operations per second.  The column ``$\times 1$'' indicates operations performed singly and ``$\times 100$'' indicates operations performed in batches of 100.  These figures are for random additions; doubling is $\approx$ 5\% slower in genus 2 and $\approx$ 2\% faster in Genus 3.  All values are for single-threaded execution.\\

The columns $E$ and $S$ show the throughput of the exponentiation ($E$) and primorial-steps search ($S$) performed by Algorithm \ref{algorithm:Exponent}.  Values were obtained by dividing the elapsed time of a single thread (including all overhead) by the number of group operations.
\normalsize
\end{minipage}
\end{table}

\begin{table}
\begin{center}
\begin{tabular}{@{}rrrrrrr@{}}
\toprule
$n$ & $u$ & $\Pr[A]$ & $\Pr[B_{2/1}|A]$ & $\Pr[B_{4/2}|A]$ &  $\Pr[B_{3/1}|A]$ &  $\Pr[\tilde{B}_{3/1}|A]$\\
\midrule
48 & 2.0 & 100 	& 4.4 & 3.8 & 5.4 & 5.4\\
48 & 3.0 & 48    & 3.8 & 3.6 & 5.1 & 5.3\\
48 & 4.0 & 11    & 3.1 & 3.4 & 4.8 & 5.1\\
%48 & 5.0 & 1.6   & 2.4 & 3.1 & 4.2 & 4.9\\
100 & 5.4 & 0.57 & 4.1 & 4.1 & 5.2 & 5.1\\
100 & 5.4 & 0.58 & 4.1 & 4.0 & 4.8 & 5.2\\
\midrule
48 & 2.0 & 100    & 8.2 & 6.6 & 6.3 & 6.3\\
48 & 3.0 & 40     & 7.7 & 6.6 & 6.0 & 6.2\\
48 & 4.0 & 7.8    & 7.2 & 6.5 & 5.9 & 6.1\\
100 & 5.4 & 0.44	& 6.5 & 5.6 & 5.6 & 5.2\\
100 & 5.4 & 0.44	& 6.9 & 5.6 & 5.3 & 4.9\\
%48 & 5 & 0.87   & 6.5 & 6.0 & 5.4 & 6.0\\
\bottomrule
\vspace{6pt}
\end{tabular}
\caption{Jacobian Order Distributions in Genus 2 (percent)}\label{table:Distribution}
\end{center}
\begin{minipage}{1.0\linewidth}
\small
The event $A$ occurs when $\#J(C)\approx 2^n$ is $2^{n/u}$-easy (Definition \ref{definition:Easy}).  The event $B_{a/b}$ occurs when $\#J_{a/b}(C)$ contains a prime factor at least 95\% the size of $\#J_{a/b}(C)$.  For a random integer, the probability of this event is $\sim \log (20/19) \approx 5.1\%$.\\

Each row reflects a dataset of $10^6$ Jacobians, with the bottom datasets containing only Jacobians with odd order (no 2-torsion).  The datasets for $n=48$ used curves with random coefficients over a field $\Fp$ with the prime $p \in [2^{24}-2^{16},2^{24}+2^{16}]$ chosen at random.  The datasets for $n=100$ used a fixed prime $p=2^{50}-27$ with either random curve coefficients (first entry), or the parameterized family $y^2 = x^5 + 2x^3 + 7x^2 + x + t$ (second entry).
\normalsize
\end{minipage}
\end{table}

\begin{table}
\begin{center}
\begin{tabular}{@{}lrrrrrrrrr@{}}
\toprule
$n$ & $w$ & $u$ & $1/\sigma(u)$ &$B$& $E$ & $S$ & $E+S$ & $(E+S)/\sigma(u)$ & $16S$\\
\midrule
100 & 5 & 5.38 &     195 & .39 & .60  & .25 & .85 & $1.7\times 10^2$ &     4\\
110 & 5 & 5.57 &     309 & .88 & 1.3  & .57 & 1.9 & $5.9\times 10^2$ &     9\\
120 & 5 & 5.75 &     484 & 1.9 & 2.9  & 1.2 & 4.1 & $2.0\times 10^3$ &    20\\
130 & 6 & 5.92 &     745 & 4.1 & 6.3  & 2.5 & 8.7 & $6.5\times 10^3$ &    40\\
140 & 6 & 6.01 &     936 & 10 & 16  & 6.4 & 22 & $2.0\times 10^4$ &   102\\
150 & 6 & 6.25 &    1765 & 17 & 25  & 10 & 36 & $6.3\times 10^4$ &   166\\
160 & 6 & 6.40 &    2640 & 34 & 50  & 21 & 71 & $1.9\times 10^5$ &   333\\
170 & 7 & 6.55 &    3972 &   65 & 97  & 39 & 136 & $5.4\times 10^5$ &   625\\
180 & 7 & 6.70 &    6012 &  122 & 183  & 74 & 256 & $1.5\times 10^6$ &  1176\\
190 & 7 & 6.84 &    8897 &  230 & 344  & 138 & 482 & $4.3\times 10^6$ &  2212\\
200 & 7 & 7.01 &   14355 &  388 & 579  & 233 & 812 & $1.2\times 10^7$ &  3728\\
\midrule
180 & 6 & 7.01 & 14355 &  54 &  80  &  33 & 114 & $1.6\times 10^6$ &   532\\
190 & 7 & 7.14 &   20943 &  102 & 153  & 62 & 215 & $4.5\times 10^6$ &   985\\
200 & 7 & 7.27 &   30553 &  191 & 286  & 115 & 401 & $1.2\times 10^7$ &  1838\\
\bottomrule
\vspace{6pt}
\end{tabular}
\caption{Search Parameter Estimates (millions)}\label{table:Semismooth}
\end{center}
\begin{minipage}{1.0\linewidth}
\small
The value $n$ is an estimate of $\lg\#J(C)$, $w$ indicates the primorial $P_w$ in Algorithm \ref{algorithm:PrimorialSteps}, and the parameter $u$ minimizes $(E+S)/\sigma(u)$, where $\sigma(u) = G(1/u,2/u)$ estimates the probability that a random integer $N\approx 2^n$ is $B$-easy.  The value $B = 2^{n/u}$ is listed in millions, as are the remaining five columns.\\

The values $E = B/\log 2$  and $S = \sqrt{2P_w/\varphi(P_w)}B$ estimate the group operations required for exponentiation ($E$) and a primorial-steps search ($S$) with parameter $B$.  Their sum approximates the group operations used by Algorithm 1 in an unsuccessful attempt to compute $\#J(C)$, and $(E+S)/\sigma(u)$ is a heuristic estimate of the average number of group operations required per successful computation of $\#J(C)$.\\

The last column estimates the memory used by Algorithm \ref{algorithm:PrimorialSteps} assuming both sets $\mathcal{B}$ and $\mathcal{G}$ are explicitly computed, using 64 bits per tuple (see Section \ref{section:Implementation}).  The last three rows show the impact of increasing $u$ to reduce space.  If only $\mathcal{B}$ is stored, the figures in the last column should be divided by 2.
\normalsize
\end{minipage}
\end{table}

\begin{table}
\begin{center}
\begin{tabular}{@{}rrrrrrrrrr@{}}
\toprule
Genus & $n$ & $u$ & $1/r$ && $E+S$ && $(E+S)/r$&\\
\midrule
2 & 100 & 5.38 &  172 &(-12\%) & .87 &(+2\%) & $1.5 \times 10^2$& (-10\%)\\
2 & 122 & 5.80 &  483 &(-12\%) & 4.7 &(+2\%) & $2.3\times 10^3$& (-11\%)\\
3 & 122 & 5.80 &  435 &(-21\%) & 4.8 &(+3\%) & $2.1 \times 10^3$& (-18\%)\\
3 & 150 & 6.25 & 1587 &(-20\%) & 36 &(+1\%) & $5.7\times 10^3$ &(-8\%)\\
2 & 168 & 6.50 & 3448 &(-0\%) & 129 &(+2\%) & $4.5\times 10^3$ &(+2\%)\\
2 & 178 & 6.70 & 5263 &(-12\%) & 213 &(+2\%) & $1.1\times 10^6$ &(-10\%)\\
\midrule
\bottomrule
\vspace{6pt}
\end{tabular}
\caption{Actuals vs. Estimates (millions)}\label{table:Actuals}
\end{center}
\begin{minipage}{1.0\linewidth}
\small
The value $r$ is the actual success rate achieved.  The first three tests used $10^6$ curves, while the last three used $10^5$.  Deviations from estimates are shown in parentheses.
%The values for $E+S$ and $(E+S)/r$ are in millions of group operations.
\normalsize
\end{minipage}
\end{table}

\begin{table}
\begin{center}
\begin{tabular}{@{}l@{}}
\toprule
%$C:\thickspace y^2 = x^5 + 2050668744648879665x^3 + 1948613779828075789x^2 \qquad\thickspace p = 2^{61}-1$\\
%$\qquad\qquad\quad + 777475836699358935x + 1981141960889113537$\\
%$\qquad a_1 = -602512302,\quad a_2 = -1028770098952312018$\\
%$\#J_{3/1}(C)$ and $\#J_{3/1}(\Ctwist)$ are 244-bit primes, $\#J_2(C)$ divisible by 3.\\
%\midrule
$C:\thickspace y^2 = x^5 + x + 456579, \qquad\qquad\qquad\qquad\qquad\qquad\qquad\qquad\qquad\quad p = 2^{61}-1$\\
$\qquad a_1 = 867588246,\quad a_2 = 503655589160075568$\\
$\#J_{3/1}(C)$ and $\#J_{3/1}(\Ctwist)$ are 244-bit primes, $\#J_2(C)$ not divisible by 3.\\
\midrule
$C:\thickspace y^2 = x^5 + x + 127861, \qquad\qquad\qquad\qquad\qquad\qquad\qquad\qquad\qquad\thickspace\thickspace p = 2^{84}-35$\\
$\qquad a_1 = -2092369310828,\quad a_2 = 35830907425009491385101310$\\
$\#J(C) =  2^5\cdot 3^2\cdot 1299112566516217620665269205633002367450315129777$\\
\midrule
$C:\thickspace y^2 = x^5 + x + 89993, \qquad\qquad\qquad\qquad\qquad\qquad\qquad\qquad\qquad\quad p = 2^{84}-35$\\
$\qquad a_1 = 1236014582768,\quad a_2 =  -20956811918028115290034218$\\
$\#J_{3/1}(C)$ is a 336-bit prime, $\#J_2(C)$ not divisible by 3.\\
\midrule
$C:\thickspace y^2 = x^5 + x + 202214, \qquad\qquad\qquad\qquad\qquad\qquad\qquad\qquad\qquad\quad p = 2^{89}-1$\\
$\qquad a_1 = -52033004229306,\quad a_2 =  1618004552234213280766854490$\\
$\#J(C) = 2^2\cdot 3^2\cdot 5\cdot 2128466028980222265110760419187916380742710181533203$\\
\midrule
$C:\thickspace y^2 = x^5 + x + 207686, \qquad\qquad\qquad\qquad\qquad\qquad\qquad\qquad\qquad\quad p = 2^{89}-1$\\
$\qquad a_1 = 37333142265075,\quad a_2 =  1342175488412716989278850463$\\
$\#J(\Ctwist_2) = 13^2\cdot P_{349}$ where $P_{349}$ is a 349-bit prime.\\
\midrule
$C:\thickspace y^2 = x^5 + 81x + 15466464,\qquad\qquad\qquad\qquad\qquad\qquad\qquad\qquad\quad p = 2^{89}-1$\\
$\qquad a_1 = -29105979141185,\quad a_2 =  216189507687913446441772723$\\
$\#J_{3/1}(C) = 7\cdot P_{354}$, where $P_{354}$ is a 354-bit prime, $\#J_2(C)$ not divisible by 3.\\
\midrule
$C:\thickspace y^2 = x^5 + 2x^3+3x^2+5x+1050,\qquad\qquad\qquad\qquad\qquad\qquad\quad\thickspace\thickspace p = 2^{93}-25$\\
$\qquad a_1 = 20868893099084,\quad a_2 =  14008940235908131442826126566 $\\
$\#J_{3/1}(C) = 7\cdot 313\cdot P_{361}$, where $P_{361}$ is a 361-bit prime, $\#J_2(C)$ divisible by 3.\\
\bottomrule
\vspace{6pt}
\end{tabular}
\caption{Genus 2 Examples}\label{table:Genus2Curves}
\end{center}
\begin{minipage}{1.0\linewidth}
\small
%The values of the primes not listed can be recovered from $J_{3/1}(C) = P(\omega)P(\omega^2)$, where $\omega = e^{2\pi i/3}$, and $J_2(\Ctwist_2) = P(i)P(-i)$.
The curve $\Ctwist_2$ is the quadratic twist of $C$ over $\Fptwo$.
\normalsize
\end{minipage}
\end{table}

\begin{table}
\begin{center}
\begin{tabular}{@{}l@{}}
\toprule
$C:\thickspace y^2 = x^7 + 3x^5 + x^4 + 4x^3 + x^2 + 5x + 851385,\qquad\qquad\qquad\qquad p = 2^{50}-27$\\
$\qquad a_1 = 13792821,\quad a_2 = 98748931364073,$\\
$\qquad a_3 = -4912096020329124903571$\\
$\#J(C) = 1427247710190335132030763894493884791800228867$\\
\midrule
$C:\thickspace y^2 = x^7 + 28x^3 + 18x^2 + 27x + 69621,\qquad\qquad\qquad\qquad\qquad p = 3\cdot 10^{16}+29$\\
$\qquad a_1 = -200710015,\quad a_2 = 49691549823351179,$\\
$\qquad a_3 = -9387711520293250802133155$\\
$\#J(C) = 5^2\cdot 373\cdot 2895442339877862336809237112865944284512053683$\\
\midrule
$C: y^2 = x^7 + 3x^5 + x^4 + 4x^3 + x^2 + 5x + 84538,\qquad\qquad\qquad\qquad\quad\thickspace p = 2^{61}-1$\\
$\qquad a_1 = -255251897,\quad a_2 = 3731171990845206887,$\\
$\qquad a_3 = -1915761422452218541377951998$\\
$\#J(C) =  2^4\cdot 3^5\cdot 17\cdot 223\cdot 831781325652289358544190241299568732364985371373$\\
\bottomrule
\vspace{6pt}
\end{tabular}
\caption{Genus 3 Examples}\label{table:Genus3Curves}
\end{center}
\begin{minipage}{1.0\linewidth}
\small
See \url{http://math.mit.edu/~drew/} for additional examples and verification details.
%The $L$-polynomial is given by $P(z) = p^3z^6+a_1p^2z^5+a_2pz^4+a_3z^3+a_2z^2+a_1z+1$.  All curves are defined over $\Fp$. 
\normalsize
\end{minipage}
\end{table}

%    Text of article.

%    Bibliographies can be prepared with BibTeX using amsplain,
%    amsalpha, or (for "historical" overviews) natbib style.
\bibliographystyle{amsplain}
%\bibliography{general}
\providecommand{\bysame}{\leavevmode\hbox to3em{\hrulefill}\thinspace}
\providecommand{\MR}{\relax\ifhmode\unskip\space\fi MR }
% \MRhref is called by the amsart/book/proc definition of \MR.
\providecommand{\MRhref}[2]{%
  \href{http://www.ams.org/mathscinet-getitem?mr=#1}{#2}
}
\providecommand{\href}[2]{#2}

\end{document}